\documentclass[runningheads]{llncs}
\newcommand{\dd}{\scriptsize{\textnormal{d}}}
\newcommand{\grad}{\textnormal{grad}}

\newcommand{\Symp}{\textnormal{Sym}^+}

\newcommand{\SO}{\textnormal{SO}}
\usepackage[utf8]{inputenc}
\usepackage{graphicx}
\usepackage{xcolor}
\usepackage{overpic}
\usepackage{wrapfig,lipsum}

\usepackage{amsmath, amsthm, amssymb, amsfonts}
\usepackage{array}
\usepackage{nicefrac}
\usepackage{hyperref}

\DeclareRobustCommand{\ShowColormap}{\raisebox{-0.14em}{\includegraphics[height=.8em]{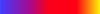}}}

\begin{document}
%
\title{Nonlinear Regression on Manifolds for\\Shape Analysis using Intrinsic Bézier Splines}
\titlerunning{Intrinsic Bézier Splines for Regression on Riemannian Manifolds}
%
\author{Martin Hanik\inst{1} \and
Hans-Christian Hege\inst{1} \and
Anja Hennemuth\inst{2} \and
Christoph von Tycowicz\inst{1}}

\authorrunning{M. Hanik et al.}
%
\institute{Zuse Institute Berlin, Berlin, Germany \\
\email{\{hanik,hege,vontycowicz\}@zib.de} \and
Institute for Imaging Science and Computational Modelling in Cardiovascular Medicine, Charité – Universitätsmedizin Berlin, Berlin, Germany 
\email{anja.hennemuth@charite.de}}
%
\maketitle              
%

\begin{abstract}
	Intrinsic and parametric regression models are of high interest for the statistical analysis of manifold-valued data such as images and shapes. The standard linear ansatz has been generalized to geodesic regression on manifolds making it possible to analyze dependencies of random variables that spread along generalized straight lines. Nevertheless, in some scenarios, the evolution of the data cannot be modeled adequately by a geodesic.
    We present a framework for nonlinear regression on manifolds by considering Riemannian splines, whose segments are Bézier curves, as trajectories.
    Unlike variational formulations that require time-discretization, we take a constructive approach that provides efficient and exact evaluation by virtue of the generalized de Casteljau algorithm.
    We validate our method in experiments on the reconstruction of periodic motion of the mitral valve as well as the analysis of femoral shape changes during the course of osteoarthritis, endorsing Bézier spline regression as an effective and flexible tool for manifold-valued regression.

\end{abstract}

\section{Introduction}
    Manifold-valued data arises in many medical applications, for example as image data or in the form of 2D/3D shapes, and sophisticated tools for its analysis have become increasingly important. Regression methods are central to modern statistics and research for their applicability to nonlinear spaces is fuelled by an ever-growing number of large longitudinal studies~\cite{gerig_longitudinal_modeling}. Consequently, geodesic regression~\cite{Fletcher13_Geodesic_Regression,Niet_geodesic_regression} was introduced as a generalization of linear regression. It allows to test whether given instances in a Riemannian manifold can be well approximated by a generalized straight line. Nevertheless, there are processes that cannot be accurately described by a geodesic, e.g., periodic motion or processes with saturation which slow down after some time. 
    In order to handle these cases, both non-parametric~\cite{Davis_kernel_regression,Mallasto_wraped_gaussian,smoothing_splines} and parametric models have been studied.
    In the latter category, Riemannian polynomials~\cite{Hinkle_ea14_Riemannian_Polynomials} and splines~\cite{SINGH_variational_spline} have been considered for nonlinear regression. They are defined, for example by employing variational principles, as solutions to differential equations involving curvature terms. Therefore, evaluation and optimization is complicated and numerically expensive since there are no closed-form solutions available in general.
    
    As an alternative, we propose to use manifold-valued Bézier curves~\cite{Nava-Yazdani_ea_De_Casteljaus,POPIEL_Bezier}. They coincide with polynomial curves in Euclidean space, are intrinsic to the manifold (i.e., independent of a choice of coordinates) and more flexible than geodesics. In contrast to Riemannian polynomials, Bézier curves allow for explicit formulas, which enables us to evaluate them directly without time-discretization. This can improve computational speed without suffering a loss of accuracy.
    Furthermore, we can combine two such curves to a differentiable spline independently of the degrees of the Bézier segments. This is again an advantage over polynomial curves.
    While variational spline models allow for piecewise composition, there is no clear way to define them for even degrees~\cite{Hinkle_ea14_Riemannian_Polynomials}. Therefore, the introduction of flexible, intrinsic splines is a key contribution of this work. Our model features closed-form, numerically stable and efficient expressions for the gradient of the regression objective in terms of concatenated adjoint Jacobi fields~\cite{Bergman_ea_Variation_Model}. In particular, we derive an algorithm that only requires basic Riemannian operations: the exponential and logarithmic map as well as certain Jacobi fields.
    Notably, closed-form expressions for these operations are available for many manifolds; in particular, they are known for Kendall's shape space~\cite{Nava-Yazdani_ea18_Geodesic_Analysis} and shape models based on differential~\cite{vTycowicz_ea_Differential_Coordinates} and fundamental~\cite{Ambellan_Surface_theoretic} coordinates.
    
    While the method can be generally applied to data on any manifold, we provide two specific examples from shape analysis.
    First, we regress the data of 100 highly resolved femur geometries with different severeness of osteoarthritis against their grade in the Kellgren Lawrence grading system. Second, we reconstruct the full motion cycle of a mitral valve from 3D geometries that were derived from ultrasound images. To the best of our knowledge, we are the first to present intrinsic regression results of such a periodic process.

\section{Spline Regression}

\subsubsection{Tools from Riemannian Geometry.}
    Before we introduce Bézier curves on manifolds we recall some important facts from Riemannian geometry; for more information see for example \cite{doCarmo}. As is often done, we use ``smooth'' synonymously with ``infinitely often differentiable''.
    
    A Riemannian manifold is a differentiable manifold $M$ together with a Riemannian metric $\langle \cdot,\cdot \rangle_p$ that assigns to each tangent space $T_p M$ a smoothly varying scalar product. As a result, a distance function $d$ is induced 
    on $M$.
    Every Riemannian manifold comes with a unique connection $\nabla$ called Levi-Civita connection. Given two vector fields $X,Y$ on $M$ it yields a natural way to differentiate $Y$ along $X$; we denote the resulting vector field by $\nabla_X Y$.

    A geodesic $\gamma$ is a generalized straight line and its defining property is vanishing of acceleration, i.e., $\nabla_{\gamma'}\gamma' = 0$, where $\gamma ' := \frac{\dd}{\dd t} \gamma$.
    An important fact is that every point in $M$ has a so-called convex neighbourhood $U$. Each pair $p,q \in U$ can be joined by a unique length-minimizing geodesic $[0,1] \ni t \mapsto \gamma(t;p,q)$ that lies completely in $U$. In the following, we always assume to work in a convex neighbourhood. Then, $\gamma$ is also differentiable with respect to its starting and end point. Explicit formula of these differentials involve the Riemannian curvature tensor $R$ (which intuitively measures local deviation from flat space; see \cite[Ch. 4]{doCarmo}). It determines Jacobi fields $J$ along $\gamma$ as solutions to the linear second order differential equation
        $\nabla_{\gamma'} \nabla_{\gamma'}J + R(J,\gamma')\gamma' = 0.$
    Considering the boundary value problem $J(0) = X$, $J(1) = 0$, we denote its solution by $J_X$.
    Then, the derivative of $\gamma$ w.r.t. its starting point $p$ in direction $X \in T_p M$ is given by $J_X$, i.e., $d_p\gamma(t;\cdot,q)(X) = J_X(t)$ for all $t\in [0,1]$.
    Furthermore, since $\gamma(t;p,q) = \gamma(1-t;q,p)$, endpoint variations are given analogously~\cite[Sec. 3.1]{Bergman_ea_Variation_Model}.

    Another important map is the Riemannian exponential. Let $X \in T_pM$ such that there is a geodesic $[0,1] \ni t \mapsto \gamma(t;p,q)$ in $U$ with $X = \gamma'(0;p,q)$. The exponential map at $p$ is then defined by $\exp_p(X) := q$. Its inverse is the Riemannian logarithm $\log_p$. In particular, we have $\log_p(q) = \gamma'(0;p,q)$. 

    The adjoint $A^*$ of a linear operator $A:T_pM \to T_qM$ is given, as usual, by the linear operator from $T_qM$ to $T_pM$ that conserves the scalar product, i.e., $\langle AX,Y \rangle_q = \langle X,A^*Y \rangle_p$ for all $X \in T_pM,\ Y \in T_qM$. 
    
    Later, we want to calculate the gradient of a composition of functions. If $f:M \to M$ and $g:M \to \mathbb{R}$ are smooth, then the chain rule for gradients reads $\grad_p (g \circ f) = {d_p f}^*(\grad_{f(p)} g)$, i.e., the gradient of $g$ at $f(p)$ is ``transported'' to the tangent space at $p$ by the adjoint differential of $f$.

\subsubsection{Bézier curves.} In the following we restrict the domain of definition to $[0,1]$ for clarity. This does not influence generality as reparametrizations are always possible. In particular, geodesics $\gamma$ can be defined on arbitrary intervals by changing the speed of travel, i.e., the length of the velocity vector.
    \begin{figure}
        \centering
        \begin{minipage}[b]{0.45\linewidth}
            \centering
            \includegraphics[width=\textwidth]{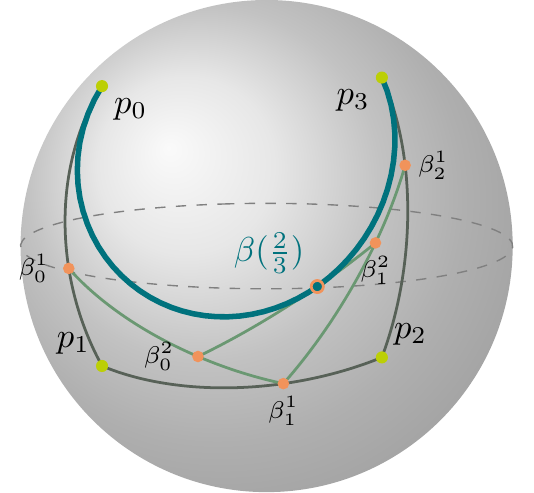}
            \vspace{.1mm}
            \end{minipage}
            \hspace{0.5cm}
            \begin{minipage}[b]{0.45\linewidth}
            \centering
            \includegraphics[width=\textwidth]{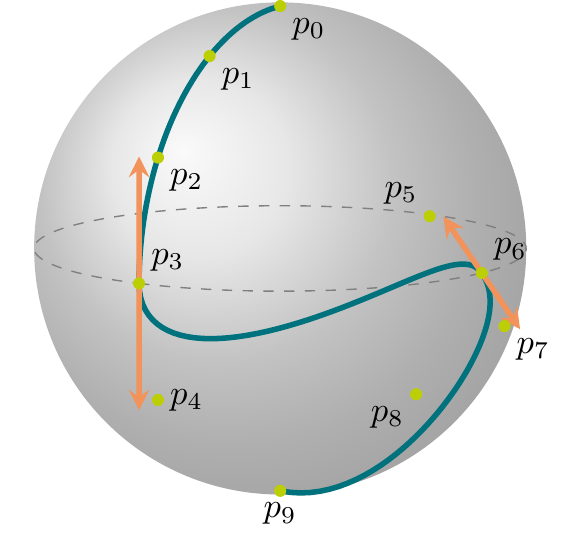}
        \end{minipage}

        \caption{Left: A cubic Bézier curve $\beta$ on the sphere $\mathcal{S}^2$ and the construction of $\beta(2/3)$ by the de Casteljau algorithm. Right: A Bézier spline with 3 cubic segments on $\mathcal{S}^2$.}
        \label{fig:Bezier_spline}
    \end{figure}
    
    A set of $k+1$ control points $p_0,\dots,p_k \in U$ defines a \textit{Bézier curve} $\beta:[0,1] \to M$ \textit{of order} $k$ according to the \textit{generalized de Casteljau algorithm}
    \begin{align}
        \beta_i^0(t) &:= p_i, \nonumber \\
        \beta_i^l(t) &:= \gamma(t;\beta_i^{l-1}(t),\beta_{i+1}^{l-1}(t)), \quad l=1,\dots,k, \quad i=0,\dots,k-l,
        \label{eq: deCasteljau}
    \end{align}
    by $\beta(t) := \beta_0^k(t).$
     
    Note that $\beta(0) = p_0$ and $\beta(1) = p_k$. Furthermore, the velocities of $\beta$ at these points are
    \begin{equation} \label{eq: velocities}
        \beta'(0) = k \log_{p_0}(p_1) \text{ and } \beta'(1) = -k \log_{p_k}(p_{k-1});
    \end{equation}
    see \cite[Thm. 1]{POPIEL_Bezier}. The algorithm is visualized on the left of Fig.~\ref{fig:Bezier_spline}. Whenever of interest, we will make the dependence of $\beta$ on its control points explicit by writing $\beta(t;p_0,\dots,p_k)$.
    Note that if there are only 2 control points $p_0,p_1$, then $\beta$ is just the geodesic from $p_0$ to $p_1$. In Euclidean space, the above algorithm is the ordinary de Casteljau algorithm (because there geodesics are straight lines) and it is a well known fact that then $\beta$ is a curve with polynomials of order at most $k$ as entries. 
    
    Property (\ref{eq: velocities}) allows us to fit Bézier curves of possibly different orders together to a differentiable spline. For $i=0,\dots,L-1$ let $p^{(i)}_0,\dots,p^{(i)}_{k_i}$ be the control points of $L$ Bézier curves such that 
    \begin{equation} \label{eq: C1_condition}
        p^{(i)}_{k_i} = p^{(i+1)}_0 \quad \text{and} \quad  \gamma \left( \frac{k_i}{k_i + k_{i+1}};p^{(i)}_{k_i-1},p^{(i+1)}_1 \right) = p^{(i+1)}_0
    \end{equation}
    for all $i=1,\dots,L-2$.
    Then, we define the \textit{Bézier spline} $B$ by
    \begin{equation} \label{eq: spline}
        B(t) := \begin{cases} \beta(t;p^{(0)}_0,\dots,p^{(0)}_{k_0}), &\quad t \in [0,1],\\ \beta(t-i;p^{(i)}_0,\dots,p^{(i)}_{k_i}), &\quad t \in (i,i+1], \quad i=1,\dots,L-1. \end{cases}
    \end{equation}
    From (\ref{eq: velocities}) it follows that $B$ is $C^1$, i.e., we can make $B$ differentiable by aligning the three control points at the connections thereby removing one degree of freedom. For more details see \cite[Sec. 2.3]{Gousenbourger_data_fitting}.  Note that we could add further restrictions to ensure that $B$ is $C^2$ \cite[p. 119]{POPIEL_Bezier}.
        
    If $L > 1$ and the first and last segment of $B$ are at least cubic, we can consider closed Bézier splines. Then, $B$ is $C^1$ and closed if and only if (\ref{eq: C1_condition}) extends cyclically, that is, we also have
    $$p^{(L-1)}_{k_{L-1}} = p^{(0)}_0 \quad \text{and} \quad  \gamma \left( \frac{k_{L-1}}{k_{L-1} + k_{0}};p^{(L-1)}_{k_{L-1}-1},p^{(0)}_1 \right) = p^{(0)}_0.$$
    
    In the following, we set 
    \begin{equation*}
    K := \begin{cases} k_0+k_1+\cdots+k_{L-2}+k_{L-1}, &\quad B \text{ non-closed,} \\
    k_0+k_1+\cdots+k_{L-2}+k_{L-1}-1, &\quad B \text{ closed,} \end{cases}
    \end{equation*}
    and denote the set of $K+1$ \textit{distinct} control points of $B$ by $p_0,\dots,p_K$. In the non-closed case this means
    $$(p_0,\dots,p_K) := \left(p^{(0)}_0,\dots,p^{(0)}_{k_0},p^{(1)}_1,\dots,p^{(1)}_{k_1},\dots,p^{(L-1)}_{1},\dots,p^{(L-1)}_{k_{L-1}} \right) \in M^{K+1},$$ 
    while $p_0^{(0)}$ is left out for closed $B$.
    An example of a $C^1$ spline with three cubic segments and 10 distinct control points is shown on the right of Fig.~\ref{fig:Bezier_spline}.

\subsubsection{The Model.}
    Let $N$ data points $q_i \in U$ with corresponding scalar parameter values $t_i$ (for example points in time) be given. We suppose that the data points $q_i$ are realizations of an $M$-valued random variable $Q$ that depends on the deterministic variable $t \in \mathbb{R}$ according to the model
    $$Q(t) = \exp_{B(t;p_0,\dots,p_K)}(\epsilon).$$
    Here, $\epsilon$ is a random variable that takes values in the tangent space $T_{B(t)}M$. The control points $p_0,\dots,p_K$ are the unknown parameters. In Euclidean space it reduces to polynomial spline regression since Bézier curves and polynomials coincide. Note that our model is a generalization of geodesic regression~\cite{Fletcher13_Geodesic_Regression}, which it reduces to when $B$ consists of a single segment with 2 control points.

\subsubsection{Least Squares Estimation.} 
    Given $N$ realizations $(t_j,q_j) \in \mathbb{R} \times U$, the \textit{sum-of-squared error} is defined by
    \begin{equation} \label{eq: error}
          \mathcal{E}(p_0,\dots,p_K) := \frac{1}{2} \sum^N_{j=1} d \Big(B(t_j;p_0,\dots,p_K),q_j \Big)^2.
    \end{equation}
    Then, we can formulate a least squares estimator of the Bézier spline model as the minimizer of this error, which under certain conditions agrees with the maximum likelihood estimation~\cite{Fletcher13_Geodesic_Regression}. 
     We would like to emphasize that none of the control points agrees with the data points as is the case in spline interpolation.
      
    In general, minimizers of~(\ref{eq: error}) are not known analytically, which makes iterative schemes necessary. Therefore, we apply Riemannian gradient descent. (For optimization on manifolds see \cite{optimazation_manifolds}.) 
    The gradient of $\mathcal{E}$ can be computed w.r.t. each control point individually.
    We write $B_t(p_i)$ for the map $p_i \mapsto B_t(p_i) := B(t;p_0,\dots,p_i,\dots,p_k)$ and define the functions $p \mapsto \tau_j(p) := d(p,q_j)^2$. It is known that $\grad_p \tau_j = -2\log_p(q_j)$ for each $p \in U$.
    When we consider the $j$-th summand on the right-hand side of (\ref{eq: error}), the chain rule implies that its gradient w.r.t. the $i$-th control point is given by 
    $$\grad_{p_i} (\tau_j \circ B_{t_j}) = d_{p_i}B_{t_j}^* \left( \grad_{B_{t_j}(p_i)} \tau_j \right) = -2\ d_{p_i}B_{t_j}^* \left(\log_{B_{t_j}(p_i)}(q_j) \right).$$
    Using (\ref{eq: error}) then gives the gradient of $\mathcal{E}$.
    The operator $d_{p_i}B_t^*$ ``mirrors'' the construction of the segment of $B$ to which $p_i$ belongs by transporting the vector $\log_{B_{t_i}(p_i)}(q_j)$ backwards along the ``tree of geodesics'' defined by the de Casteljau algorithm~(\ref{eq: deCasteljau}). More precisely, the result is a sum of vectors in $T_{p_i}M$ that are values of concatenated adjoint differentials of geodesics w.r.t.\ starting and end point. In symmetric spaces, for example, they are known in closed form. For a detailed inspection of $d_{p_i}B_t^*$ we refer to~\cite[Sec. 4]{Bergman_ea_Variation_Model}.

    As initial guess for the gradient descent, we choose $(p_0,\dots,p_K)$ along the geodesic polygon whose corners interpolate the data points that are closest to knot points w.r.t.\ time.
    
\section{Experiments}
    Although physical objects themselves are embedded in Euclidean space, their \textit{shape features} are best described by more general manifolds requiring Riemannian geometric tools for statistical analysis thereon; see for example~\cite{BauerBruverisMichor2014,Kendall_ea_Shape_Theory,vTycowicz_ea_Differential_Coordinates}.
    To test our regression method for shape analysis, we apply it to two types of 3D data: (i) distal femora and (ii) mitral valves given as triangulated surfaces. We perform the analysis in the shape space of \textit{differential coordinates}~\cite{vTycowicz_ea_Differential_Coordinates}. That is, for homogeneous objects given as triangular meshes in correspondence, we choose their intrinsic mean~\cite{fletcher_stat} as reference template and view all objects as deformations thereof. (We assume that the meshes are rigidly aligned, e.g., by generalized Procrustes alignment \cite{generalized_procrustes}.) On each face of a mesh, the corresponding deformation gradient is constant and, therefore, can be encoded as a pair of a rotation and a stretch, i.e., as an element of the Lie group of 3 by 3 rotation and symmetric positive definite matrices $\SO(3) \times \Symp(3)$. Denoting the Frobenius norm by $\|\cdot\|_F$, metrics are chosen such that the distance functions become $d_{SO}(R_1,R_2) := \| \log(R_1^TR_2) \|_F$ and $d_{\Symp}(S_1,S_2) := \|\log(S_2)-\log(S_1) \|_F$, respectively. Suppose the number of triangles per object is $m$, then the full shape space is the product space $(\SO(3) \times \Symp(3))^m$. Using the product metric, statistical analysis thereon can be done face-wise and separately for rotations and stretches. We implemented a prototype of our method in MATLAB using the MVIRT toolbox~\cite{Bergmann_MVIRT}.

\subsubsection{Distal femora.} 
   Osteoarthritis (OA) is a degenerative disease of the joints that is, i.a.,\ characterized by changes of the bone shape.
   To evaluate our model, we regress the 3D shape of distal femora against OA severity as determined by the Kellgren-Lawrence~(KL) grade~\cite{Kellgren1957RadiologicalAO}---an ordinal scale from 0 to 4 based on radiographic features.
   Our data set comprises 100 shapes (20 per grade) of randomly selected subjects from the OsteoArthritis Initiative (a longitudinal, prospective study of knee OA) for which segmentations of the respective magnetic resonance images are publicly available (\href{https://doi.org/10.12752/4.ATEZ.1.0}{https://doi.org/10.12752/4.ATEZ.1.0})~\cite{AMBELLAN_data}.
   In a supervised post-process, the quality of segmentations as well as the correspondence of the extracted triangle meshes (8,988 vertices\;/\;17,829 faces) were ensured.
   
    \begin{center}
        \begin{table}[h]
            \centering
            \begin{tabular}{ >{\centering\arraybackslash} p{4cm}  p{1.5cm}  p{1cm} }
                Order of Bézier curve & $R^2$ & $R^2_{rel}$ \\
                \hline
                1 & 0.05 & 0.57 \\ 
                2 & 0.07 & 0.78 \\ 
                3 & 0.08 & 0.90 \\ 
            \end{tabular}
            \caption{The computed $R^2_{\textnormal{rel}}$ and $R^2$ statistics of the regressed (w.r.t. KL grade) geodesic, quadratic and cubic Bézier curve for data of distal femora}
            \label{tab: results_femur}
        \end{table}
    \end{center}
  
    For $i=0,\dots,4$, the shapes with grade $i$ are associated with the value $t_i = \nicefrac{i}{4}$.
    We use our method to compute the best-fitting geodesic, quadratic and cubic Bézier curve.
    In order to compare their explanatory power, we calculate for each the corresponding manifold-valued $R^2$ statistic that, for $r_1,\dots,r_N \in M$ and \textit{total variance} $\textnormal{var}\{r_1,\dots,r_N\} := \nicefrac{1}{N} \min_{q \in M} \sum_{j=1}^N d(q,r_j)^2$, is defined by~\cite[p. 56]{fletcher_stat}
    $$R^2 = 1- \frac{\textnormal{unexplained variance}}{\textnormal{total variance}}:=1-\frac{\nicefrac{2}{N} \mathcal{E}(\beta)}{\textnormal{var}\{r_1,\dots,r_N\}} \in [0,1].$$
    The statistic measures how much of the data's total variance is explained by $\beta$.
   
   For $j=1,\dots,20$  and $l=0,\dots,4$, let $q_j^{(l)}$ be the $j$-th femur shape with KL grade $l$. 
   Note that, for the described setup, the unexplained variance is bounded from below by the sum of the per-grade variances,
   i.e., $\sum_{l=0}^4\textnormal{var}\{q_1^{(l)},\dots,q_{20}^{(l)}\}.$
   In particular, for our femur data this yields an upper bound for the $R^2$ statistic of $R^2_{\textit{opt}} \approx 0.0962$. 
   Hence, we also provide relative values $R^2_{\textit{rel}} := \nicefrac{R^2}{R^2_{\textit{opt}}}$ for comparison.
   The results are shown in Table \ref{tab: results_femur}.
   
   \begin{figure}[t]
     \begin{center}
     	\includegraphics[width=\textwidth]{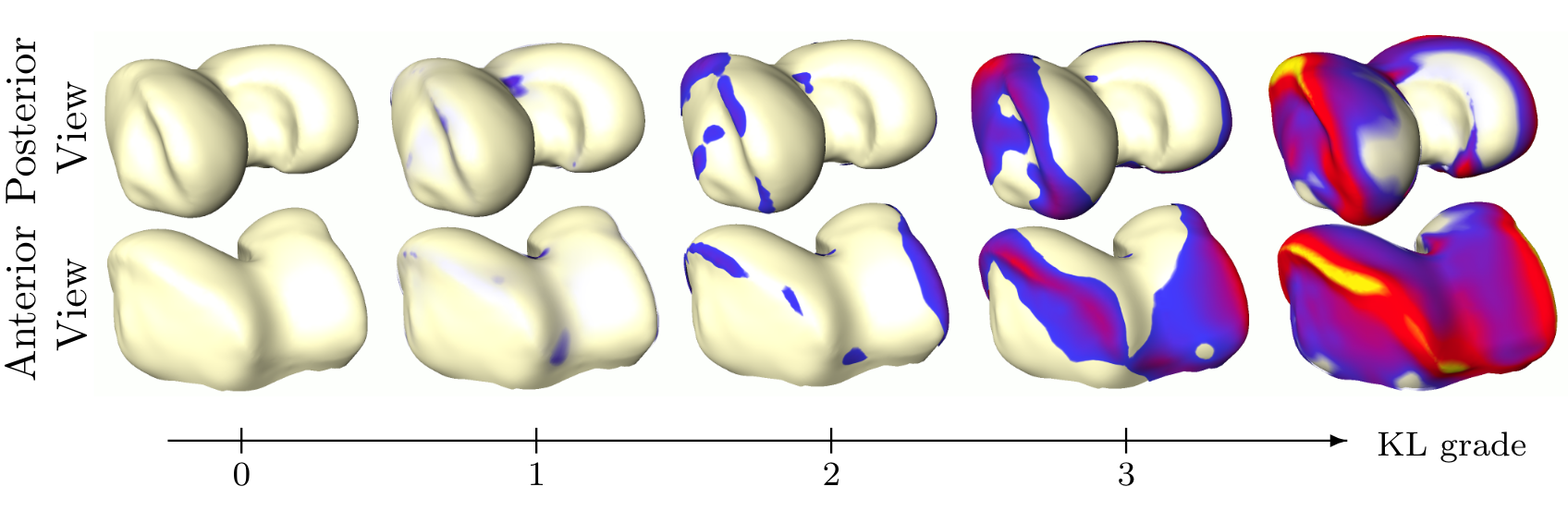}
    	\caption{Cubic regression of distal femora.
    	Healthy regressed shape ($\textnormal{KL}=0$) together with subsequent grades overlaid wherever the distance is larger than 0.5mm, colored accordingly (0.5~\ShowColormap~3.0).
    	}
     \label{fig:femur_cubic}
     \end{center}
    \end{figure}

    The computed cubic Bézier curve is displayed in Figure \ref{fig:femur_cubic}. The obtained shape changes consistently describe OA related malformations of the femur, viz., widening of the condyles and osteophytic growth. Furthermore, we observe only minute bone remodeling for the first half of the trajectory, while accelerated progression is clearly visible for the second half. The substantial increase in $R_{\textnormal{rel}}^2$ suggests that there are nontrivial higher order phenomena involved which are captured poorly by the geodesic model. Moreover, as time-warped geodesics are contained in the search space we can inspect time dependency. Indeed, for the cubic femoral curve the control points do not belong to a single geodesic, confirming higher order effects beyond reparametrization. 
    
\subsubsection{Mitral Valve.}

    Diseases of the mitral valve such as mitral valve insufficiency (MI) are often characterized by a specific motion pattern and the resulting shape anomalies can be observed (at least) at some point of the cardiac cycle. In patients with MI, the valve's leaflets do not close fully or prolapse into the left atrium during systole. Blood then flows back lowering the heart's efficiency. 
    
    We compute regression with Bézier splines for the longitudinal data of a diseased patient's mitral valve. Sampling the first half of the cycle (closed to fully open) at equidistant time steps, 5 meshes (1,331 vertices / 2,510 faces) were extracted from a 3D+t transesophageal echocardiography (TEE) sequence as described in~\cite{tautz_MVdata}. Let $q_1,\dots,q_5$ be the corresponding shapes in the space of differential coordinates. In order to approximate the full motion cycle we use the same 5 shapes in reversed order as data for the second half of the curve. 
    Because of the periodic behaviour, we choose a closed spline with two cubic segments as model and assume an equidistant distribution of the data points along the spline, i.e., we employ $\{(0,q_1), (\nicefrac{1}{4},q_2), (\nicefrac{1}{2},q_3), (\nicefrac{3}{4},q_4), (1,q_5),(\nicefrac{5}{4},q_4),(\nicefrac{3}{2},q_3),(\nicefrac{7}{4},q_2)\}$ as the full set. 
    
    \begin{figure}[t]
        \begin{center}
     	    \includegraphics[width=\textwidth]{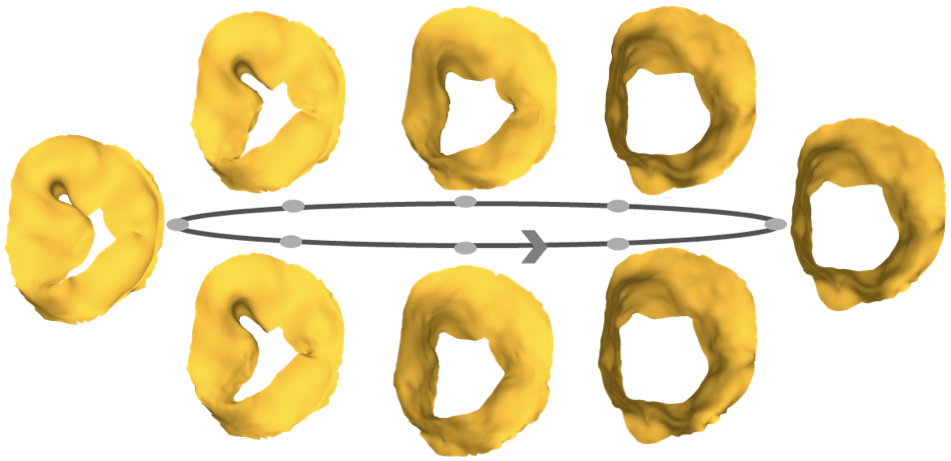}
    	    \caption{Reconstructed meshes from regression of longitudinal mitral valve data covering a full cardiac cycle. The spline consists of two cubic segments.}
            \label{fig:valve}
         \end{center}
    \end{figure}
    
    The regressed cardiac trajectory is shown in Fig.~\ref{fig:valve}.
    Our method successfully estimates the valve's cyclic motion capturing the prolapsing posterior leaflet. 
    It shows the potential for improved reconstruction of mitral valve motion in presence of image artifacts like TEE shadowing and signal dropout.
    This, in turn, facilitates quantification of geometric indices of valve function such as orifice area or tenting height.
    
    \section{Conclusion}
    
    We presented a parametric regression model that combines high flexibility with efficient and exact evaluation. In practice, it can be used for many types of manifold-valued data as it relies only on three basic differential geometric tools that can be computed explicitly in many important spaces. In particular, we have presented two applications to shape data where we could model higher order effects and cyclic motion. 
    A remaining question is, which Bézier spline (number of segments and their order) to choose for the analysis of a particular data set. This problem of model selection poses an interesting avenue for future work.
    Moreover, we plan to extend the proposed framework to a hierarchical statistical model for the analysis of longitudinal shape data, where subject-specific trends are viewed as perturbations of a population-average trajectory represented as Bézier spline.
    
\subsubsection{Acknowledgments.}  
    M. Hanik is funded by the Deutsche Forschungsgemeinschaft (DFG, German Research Foundation) under Germany’s Excellence Strategy – The Berlin Mathematics Research Center MATH+ (EXC-2046/1, project ID: 390685689). Furthermore we 
    are grateful for the open-access dataset OAI 
    \footnote{The Osteoarthritis Initiative is a public-private partnership comprised of five contracts (N01-AR-2-2258; N01-AR-2-2259; N01-AR-2-2260; N01-AR-2-2261; N01-AR-2-2262) funded by the National Institutes of Health, a branch of the Department of Health and Human Services, and conducted by the OAI Study Investigators. Private funding partners include Merck Research Laboratories; Novartis Pharmaceuticals Corporation, GlaxoSmithKline; and Pfizer, Inc. Private sector funding for the OAI is managed by the Foundation for the National Institutes of Health. This manuscript was prepared using an OAI public use data set and does not necessarily reflect the opinions or views of the OAI investigators, the NIH, or the private funding partners.} 
    and the open-source software MVIRT \cite{Bergmann_MVIRT}.
%
%
%
    \bibliographystyle{splncs04}
    \bibliography{bibliography}
    
\end{document}